\documentclass[12pt,twoside]{article}
\usepackage{amsmath}
\usepackage{amssymb}
\usepackage{enumerate,epsfig}
\newtheorem{lemma}{Lemma}[section]
\newtheorem{theorem}{Theorem}[section]
\newtheorem{corollary}[lemma]{Corollary}

\newtheorem{prop}{Proposition}[section]

\newcommand{\R}{ {\mathbb R} }

\makeatletter
\@addtoreset{equation}{section}
\makeatother

\newcommand{\cqfd}{{\unskip\kern 6pt\penalty 500
\raise -2pt\hbox{\vrule\vbox to 6pt{\hrule width 6pt
\vfill\hrule}\vrule}\par}}
\newcommand{\ep}{{\varepsilon}}

\begin{document}
\title{Compactness for nonlinear transport equations}

\author{Fethi Ben Belgacem$^1$, Pierre-Emmanuel Jabin$^{2}$}

\footnotetext[1]{Universit\'e de Monastir,  Institut sup\'erieure
  d'informatique et de math\'ematiques de Monastir, d\'epartement de
  math\'ematiques, E-mail: fethi.benbelgacem@fst.rnu.tn,
  belgacem.fethi@gmail.com }

\footnotetext[2]{TOSCA project-team, INRIA Sophia Antipolis --
  M\'editerran\'ee, 2004 route des Lucioles, B.P.\ 93, 06902 Sophia
  Antipolis Cedex, France, and 
Laboratoire J.-A. Dieudonn\'e, Universit\'e de Nice --
  Sophia Antipolis, Parc Valrose, 06108 Nice Cedex 02, France, E-mail:
  \texttt{jabin@unice.fr}}
\date{}

\maketitle

\noindent{\bf Abstract.} We prove compactness and hence existence for
solutions to a class of non linear transport equations. The
corresponding models combine the features of linear transport
equations and scalar conservation laws. We introduce a new method
which gives quantitative compactness estimates compatible with both
frameworks.

\section{Introduction}
Recent developments in the modeling of various complex transport
phenomena (from bacteria to pedestrians' flows) have produced new and
challenging equations. In particular those models have a very
different behaviour from the usual fluid dynamics when the density is
locally high, usually as a consequence of a strict bound on the
maximum number of indivivuals that one can have at a 
given point. 

The mathematical theory for well posedness and particularly existence
is still however lacunary for those equations. The aim of this article
is thus to provide a unified framework for a general class of
conservation laws, including many of these recent additions.
More precisely, we study equations 
of the form, 
\begin{equation}
\partial_t n(t,x)+\hbox{div}\, (a(t,x)\,f(n(t,x)))=0,\quad t\in
\R_+,\ x\in \R^d,\label{maineq}
\end{equation}
where $n$ usually represents a density of individuals. $f\in
W^{1,\infty}(\R,\;\R)$ is a given function which takes local, non linear effects
into account. A typical example for $f$ is the logistic $f(n)=n\,(1-n/\bar
n)_+$, which limits the velocity of individuals when
their density is too high thus ensuring that the density never exceeds
a critical value $\bar n$. The field $a: \R_+\times \R^d\rightarrow \R^d$
provides the direction for the movement of individuals. 

Depending on
the exact model, $a$ can either be given or be related to $n$. Many
such models have been introduced in the past few years in various
contexts from chemotaxis for cells and bacteria to pedestrian flow
models. We only give here a few such examples.

Typically $a$ incorporates some non local effects on the density such
as with a convolution $a=K\star n$ or a Poisson eq. 
\begin{equation}
a(t,x)=-\nabla_x \phi(t,x),\quad -\Delta_x \phi(t,x)
=g(n(t,x)),\label{poisson}
\end{equation}
where $g$ is another given function of $n$. Such a model was
introduced in two dimensions and in the context of swarming in
\cite{TB}. The same kind of models was studied in \cite{BDS} and
\cite{DP} for chemotaxis (and typically for $g(n)=n$). 

More complicated relations between $a$ and $n$ are possible, for instance a
Hamilton-Jacobi equation as in \cite{DMP}
\begin{equation}
a(t,x)=-\nabla_x \phi(t,x),\quad -\Delta_x \phi(t,x)+\alpha\,|\nabla
\phi|^2
=g(n(t,x)),\label{hj}
\end{equation}
with $\alpha\geq 0$ (possibly vanishing)
and again $g$ a given non linear function.
 
\medskip

Eq. \eqref{maineq} can be seen as a hybrid model, combining features
of usual linear transport equation and scalar conservation laws.

Let us briefly discuss the main difficulty in obtaining existence of
distributional solutions to \eqref{maineq}.
With reasonable assumptions (like $f\sim n(1-n)$), 
it is easy to show that the density $n$
is bounded in every $L^p$ spaces. However contrary to linear transport
equations, a bound on $n$ is not enough to pass to the limit in
the nonlinear term $f(n)$ (or $g(n)$ if \eqref{poisson} or \eqref{hj}
is used). 

With Eq. \eqref{poisson} or \eqref{hj} and $n\in L^1\cap L^\infty$,
one can easily get $a\in W^{1,p}$ for any $1<p<\infty$. From that one
may obtain compactness on $a$ in $L^1_{loc}$. 

Hence as a non linear model, the main difficulty in obtaining existence of
solutions to \eqref{maineq} is to prove compactness for the density
$n$. Below we briefly indicate why the usual methods for conservation
laws do not work in this setting (see \cite{Da} or \cite{Se} for more
on conservation laws).

When $a$ is regular enough (Lipschitz more precisely), then the
usual method of compactness for scalar conservation laws
work and one can for example show propagation of BV bounds on $n$. 
Unfortunately this Lipschitz bound does not hold here in
general (only $W^{1,p}$, $p<\infty$ as explained above).
Such $BV$ bounds on $n$ can in fact only be propagated for short times
(see \cite{BDS} for instance).

For scalar conservation laws, another way to obtain
compactness is
either by compensated compactness or other regularizing
effects. However in dimension larger than $1$, those
cannot be used as the flux cannot 
be genuinely non linear (it is
in only one direction, the one given by $a$). The 1-dimensional case
is quite particular (not only in this respect) and many well posedness
results have already been obtained (see for instance \cite{DMP}).

As far as we know, \cite{DP} is for the moment
the only result showing existence to an equation like
\eqref{maineq} over any time interval and any dimension. The authors
use a kinetic formulation of \eqref{maineq}, which simply generalizes
the kinetic formulation of scalar conservation laws introduced in
\cite{LPT} (see also \cite{LPT2} and \cite{Pe}). A rigidity property
inherent to the kinetic formulation then provides compactness. 
However a precise connection between $a$ and $n$ is needed; more
precisely the result is obtained only for the case of \eqref{poisson}
(with $g=Id$ though it can obviously be extended to any $g$ suitably
regular).

We conclude this brief summary of the various techniques already in
use by mentioning gradient flows. In the context of the non
linear model \eqref{maineq}, the theory is essentially still in
development. It requires a lot of structure on the equations and that
essentially means for the moment Eq. \eqref{poisson}
{\em with} $g=Id$ (any generalization to non linear $g$ would be
problematic). We refer in particular to \cite{DNS} where the right metric
for the problem and its properties are introduced and studied. 

Gradient flows techniques were also used in \cite{MRS} for a related
problem. In that case the corresponding transport is linear but
associated with a constraint on the maximal density. In the framework
of \eqref{maineq} that would correspond to $f(\xi)=\xi\,{\mathbb I}_{\xi<1}$. 
 

\smallskip

Let us now formulate the main results of the paper. 
Consider a vanishing viscosity approximation
\begin{equation}
\begin{split}
&\partial_t n_\ep(t,x)+\hbox{div}\,
(a_\ep(t,x)\,f(n_\ep(t,x)))-\ep^2\,\Delta_xn_\ep=0,\quad t\in
\R_+,\ x\in \R^d,\\
&n_\ep(t=0,x)=n_\ep^0(x).
\end{split}\label{viscous}
\end{equation}
Instead of assuming a precise form or relation between $a_\ep$ and
$n_\ep$, we make very general assumptions on $a_\ep$. Assume that on
$[0,\ T]$
\begin{equation}
 \exists\;p>1,\qquad 
\sup_\ep\sup_{t\in[0,\ T]} \|a_\ep(t,.)\|_{W^{1,p}(\R^d)}<\infty,\label{w11}
\end{equation}
\begin{equation}
 \sup_\ep \|\hbox{div}_x a_\ep\|_{L^\infty([0,\ T]\times
    \R^d)}<\infty\label{divb}.
\end{equation}
As for linear transport equation, an additional condition is needed on
the divergence to obtain compactness. In order to be compatible with
\eqref{poisson} or \eqref{hj}, we assume
\begin{equation}
\left\{\begin{aligned}
&\hbox{div}_x a_\ep=d_\ep+r_\ep\quad\hbox{with $d_\ep$ compact and }\\
&  
\exists
  C>0,\ s.t.\ \forall \ep>0, \forall x,y,\\
&
|r_\ep(x)-r_\ep(y)|
\leq
C\,|n_\ep(t,x)-n_\ep(t,y)|.
\end{aligned}\right.\label{divcomp}
\end{equation}
Then one can prove
\begin{theorem}
Assume \eqref{w11}, \eqref{divb}, \eqref{divcomp}, that $a_\ep$ is
compact in $L^p$, that $n^0_\ep$ is
uniformly bounded in $L^1\cap L^\infty(\R^d)$ and is compact in
$L^1(\R^d)$.
Then the solution $n_\ep(t,x)$ to \eqref{viscous}
is compact in $L^1_{loc}([0,\ T]\times\R^d)$. \label{mainth}
\end{theorem}

This in particular implies existence results like
\begin{corollary}
Assume that $f\in W^{1,\infty}$, $g\in C^2$, $f(0)=g(0)=0$ 
and that $f(\xi)\,g(\xi)\geq
-C\,|\xi|$ for some given constant $C$. 
Let $n^0\in L^1\cap L^\infty(\R^d)$, $\alpha\geq 0$, then $\exists n\in
L^\infty([0,\ T],\ L^1\cap L^\infty(\R^d))$ solution in the sense of
distribution to \eqref{maineq} with \eqref{poisson}. Moreover $n$ is an
entropy solution to \eqref{maineq} in the usual sense that 
$\forall \phi\in C^2$ convex, $\exists q\in C^1$ s.t.
\[
\partial_t (\phi(n(t,x)))+{\rm
  div}_x\,(a(t,x)\,q(n(t,x)))+(\phi'(n)\,f(n)-q(n))\,{\rm div}_x\, a
\leq 0.
\]\label{existence}
\end{corollary}
Note that this is just one example of possible results, it can for
instance easily
be generalized to \eqref{hj} under corresponding assumptions. Once $a$
is given and in $W^{1,p}$ the uniqueness of the entropy solution to
\eqref{maineq} is actually not very difficult. However uniqueness
for a coupled system like \eqref{maineq}-\eqref{hj} is more delicate
and left open here.

To prove Th. \ref{mainth}, 
 we develop a new method which is a sort of quantified
version of the theory of renormalized solution and compatible with the
usual $L^1$ contractivity argument for scalar laws.

Renormalized solutions were introduced in \cite{DL} to prove
uniqueness to solutions of linear transport equations
\[
\partial_t n+{\rm div}\,(a\,n)=0.
\]
The compactness of a sequence of bounded solutions is obtained as a
consequence of the uniqueness (by proving for instance that $w-\lim_k
n_k^2=(w-\lim_k n_k)^2$). The theory was developped in \cite{DL} for $a\in
W^{1,1}$ with ${\rm div}\,a\in L^\infty$. It was later extended to
$a\in BV$, first for the particular case of kinetic equations in
\cite{Bo} (see also \cite{CJ} for the kinetic case with less than one
derivative on $a$). The general case was dealt with 
in \cite{Am} (see also
\cite{CL}). For more about renormalized solutions we refer to
\cite{ADM} and \cite{DeL}.   

The usual proof of the renormalization property relies on a commutator
estimate. It is this estimate that we have to quantify somehow here.
More precisely we try to bound quantities like
\begin{equation}
\|n_\ep\|_{p,h}^p=\int_{\R^{2d}} \frac{{\mathbb I}_{|x-y|\leq
    1}}{(|x-y|+h)^d}\,|n_\ep(t,x)-n_\ep(t,y)|^p\,dx\,dy,
\label{defnorm}\end{equation}
uniformly in $h$. Those norms can be seen as a generalization of usual
Sobolev norm, in particular we recall that 
\[
\int_{\R^{2d}} \frac{{\mathbb I}_{|x-y|\leq
    1}}{|x-y|^{d+2s}}\,|n_\ep(t,x)-n_\ep(t,y)|^2\,dx\,dy
\]
is equivalent to the usual $\dot H_s$ norm for $s\in ]0,\ 1[$. This is
  wrong though for $s=0$, {\em i.e.} $\|.\|_{2,0}$ is actually
  stronger than $L^2$. In this case $p=2$, it is in fact easy to see
  in Fourier that $\|.\|_{2,0}$ more or less controls the $\log$ of a
  derivative and thus provides compactness. 

We can prove explicit estimates for the norms \eqref{defnorm}
\begin{theorem}
Assume \eqref{w11}, \eqref{divb}, \eqref{divcomp}, that $n^0_\ep$ is
uniformly bounded in $L^1\cap L^\infty(\R^d)$ and is compact in
$L^1(\R^d)$. $\exists C>0$ only depending on the uniform bounds in
$\ep$ s.t. the solution $n_\ep(t,x)$ to \eqref{viscous} satisfies for
any $t\leq T$ 
\[\begin{split}
&\int_{\R^{2d}} \frac{{\mathbb I}_{|x-y|\leq
    1}}{(|x-y|+h)^d}\,|n_\ep(t,x)-n_\ep(t,y)|\,dx\,dy
\\
&\ \leq e^{C\,t}\;\Bigg\{\int_{\R^{2d}} \frac{{\mathbb I}_{|x-y|\leq
    2}}{(|x-y|+h)^d}\,|n_\ep^0(x)-n_\ep^0(y)|\,dx\,dy\\
&\quad \qquad \quad+\int_0^t\int_{\R^{2d}} \frac{{\mathbb I}_{|x-y|\leq
    1}}{(|x-y|+h)^d}\,|d_\ep(s,x)-d_\ep(s,y)|\,dx\,dy\,ds\\
&\quad \qquad\quad+C\,\frac{\ep^2}{h^2}+C\,|\log h|^{1/\bar p}\Bigg\},
\end{split}
\]
where $\bar p=\min(2,p)$ and $1/p^*+1/p=1$.
\label{quantth}\end{theorem}
\noindent{\bf Remarks.}\\
1. Lemma \ref{compact} below shows that Theorem \ref{quantth} in fact implies
Theorem \ref{mainth} but its proof is of course more complicated.\\
2. In addition of providing an explicit rate, Theorem \ref{quantth}
does not require the compactness of the sequence $a_\ep$. Of course as
it is uniformly in $L^\infty_t W^{1,p}_x$, it is always compact in
space but not necessarily in time.\\
3. It is possible to replace \eqref{w11} by 
\[
\sup_\ep\int_0^T \|a_\ep(t,.)\|_{W^{1,p}(\R^d)}<\infty.
\]  
The estimate then uses the exponential of this quantity instead of
$e^{Ct}$.\\
4. If the sequence $\nabla a_\ep$ is equiintegrable then some kind of
rate can also be obtained. \\
5. Assumption \eqref{divcomp} can also be extended by asking $r_\ep$
to satisfy only
\[
\|r_\ep\|_{h,1}\leq C\,\|n_\ep\|_{h,1}.
\]

The norms defined by \eqref{defnorm} are in fact critical for the
problem \eqref{maineq}. Indeed \eqref{maineq} contains the case of the
linear transport equation (take $f=Id$). In this last case, one may
use the characteristics and it was
proved in \cite{CL} that one indeed propagates a sort of $\log$ of
derivative on them. If $n^0\in W^{1,p}$ then this implies a result
like Th. \ref{quantth}. Moreover at the level of the characteristics,
it is not complicated to obtain examples showing that this logarithmic gain 
is the best one can hope for.

Note that contrary to \cite{CL}, we work here at the level of the PDE;
because of the shocks, the characteristics cannot be used when $f$ is
non linear. This unfortunately makes the corresponding proof considerably more
complicated and in particular it forces us to carefully track every
cancellation in the commutator estimate;  we also refer to \cite{BC} for
an example in a different linear situation where a problem of similar
nature is found. 

Th. \ref{quantth} gives a rate in $|\log h|^{1/\bar p}$ which is
probably not optimal. In the linear case $f=Id$, \cite{CD} shows that
the optimal rate is $1$. In our non linear situation, it seems
reasonable to conjecture that it should be the same (at least for
$p\geq 2$) but it is obviously a difficult question.

The proof of Th. \ref{quantth} requires the use of multilinear
singular integrals. This has been an important field of study in
itself (we quote only some results below) but
quite a few open questions remain, making the optimality of
Th. \ref{quantth} unclear. 

The first contributions for multilinear singular integrals were
essentially in dimension 1, see \cite{Ca}, \cite{CM} or
\cite{CM2}. The theory was later developed for instance in \cite{GT},
\cite{KS}, \cite{LT}. In
dimension 1, an almost complete answer was finally given in
\cite{MTT}. In higher dimension, the most complete result that we know
of, \cite{MTT2}, unfortunately does not contain the case that we have
to deal with here.

Let us conclude this introduction by mentioning two important and still open
problems. Of course many technical issues are still unresolved: The
optimal rate, the case where $a_\ep\in BV$ instead of at least
$W^{1.1}$... 

First of all, in many situations a bound on the divergence of $a_\ep$
is not available. However when $f$ is a logistic function for example,
Eq. \eqref{maineq} still controls the maximal compression, contrary to
a linear transport equation. It means that this case should actually
be easier to handle in the non linear setting.

Second some models do not provide any additional derivative on the
velocity field 
$a$. For instance in porous media, one finds the classical coupling
\[
a=-\nabla \phi,\quad {\rm div}_x (\alpha(n)\,\nabla \phi)=g,
\]
but one could also consider the non viscous equivalent of \eqref{hj}.
Of course the method presented here fails in those cases...

The next section gives a quick proof of Corollary \ref{existence}. The
next section is devoted to Th. \ref{mainth} and the last one to
Th. \ref{quantth}.

In the rest of the paper, 
$C$ will denote a generic constant, which may depend
on the time interval $[0,\ T]$ 
considered, uniform bounds on the initial data $n^0_\ep$ or on $a_\ep$
but which never depends on $\ep$ or the parameter $h$ that we will introduce.
\section{Proof of Corollary \ref{existence}}
%
Define a sequence of approximations $n_\ep$, $a_\ep$ where $n_\ep$
solves \eqref{viscous} with initial data $n^0$ 
and $a_\ep$ is obtained through $n_\ep$ by
solving \eqref{hj}.

As \eqref{viscous} is conservative then one obviously has
\[
\|n_\ep(t,.)\|_{L^1(\R^d)}=\|n^0\|_{L^1}.
\]
By the maximum principle
\[
\frac{d}{dt}\|n_\ep(t,.)\|_{L^\infty(\R^d)}\leq \|(f(n_\ep)\,{\rm div}
a_\ep)_{-}\|_{L^\infty}, 
\]
where $(.)_{-}$ denotes the negative part. Using \eqref{poisson} implies that
\[
\frac{d}{dt}\|n_\ep(t,.)\|_{L^\infty(\R^d)}\leq
\|(fg(n_\ep))_{-}\|_{L^\infty}\leq C\, \|n_\ep(t,.)\|_{L^\infty(\R^d)},
\]
by the assumption in Corollary \ref{existence}. Hence by Gronwall's
lemma, the sequence
$n_\ep$ is uniformly bounded in $L^\infty([0,\ T],\ L^1\cap
L^\infty(\R^d))$ for any $T>0$.

Thanks to $g(0)=0$, 
the usual estimate for \eqref{poisson} then gives that $a_\ep$ is
uniformly in $L^\infty([0,\ T],\ \dot W^{1,p}(\R^d))$ for any
$1<p<\infty$. \eqref{w11}, \eqref{divb}, \eqref{divcomp} are hence obviously
satisfied.

To apply Th. \ref{mainth}, it only remains to obtain the compactness of
$a_\ep$ (note that the refined Th. \ref{quantth} does not require
it). First we need an additional bound on $n_\ep$. Multiplying 
Eq. \eqref{viscous} by $n_\ep$ and integrating, one finds
\[
\ep \int_0^T \int_{\R^d} |\nabla n_\ep|^2\,dx\leq \int_{\R^d}
|n^0(x)|^2\,dx+\frac{1}{2}\int_0^T\int_{\R^d} |n_\ep|^2\,{\rm div}\,a_\ep\,dx.   
\]
Thus the previous bounds show that $\ep^{1/2}\,\nabla n_\ep$ is
uniformly bounded in $L^2$.

Now using the transport equation \eqref{viscous} and the
relation \eqref{poisson} implies for $h'=f'g'$
\[\begin{split}
\partial_t a_\ep=&\nabla\Delta^{-1} \partial_t (g(n_\ep(t,x)))=
-\nabla\Delta^{-1} {\rm div}\,(a_\ep
h(n_\ep))\\
&-\nabla\Delta^{-1}(g'f-h)(n_\ep)\,
{\rm div}\,a_\ep
+\ep\nabla g(n_\ep)-\ep \nabla\Delta^{-1}\,g''(n_\ep)\,|\nabla n_\ep|^2.
\end{split}\]
This proves
that $\partial_t a_\ep$ is uniformly bounded in $L^2([0,\ T])$ with values in
some negative Sobolev space. Therefore $a_\ep$ is locally compact in
$L^p(\R^d)$ with $p$ large enough, more precisely $p>(1-1/d)^{-1}$ by
Sobolev embeddings.

It only remains to control the behaviour at $\infty$ of $n_\ep$ and
hence $a_\ep$. By De la Vall\'ee Poussin, since $n^0\in L^1$, there
exists $\psi\in C^\infty$, convex with $\psi(x)\rightarrow +\infty$ as
$|x|\rightarrow +\infty$, $\nabla \psi\in L^\infty$ and s.t.
\[
\int_{\R^d} \psi(x)\,|n^0(x)|\,dx<\infty.
\]  
By the convexity of $\psi$, one obtains
\[
\frac{d}{dt}\int_{\R^d} \psi(x)\,|n_\ep(t,x)|\,dx\leq \int_{\R^d}
|\nabla \psi|\,|{\rm div}\, a_\ep|\,|f(n_\ep)|\,dx\leq C.
\]
This implies that $\forall t\in[0,\ T]$
\begin{equation}
\int_{\R^d} \psi(x)\,|n_\ep(t,x)|\,dx\leq C.
\label{queue}
\end{equation}
By \eqref{poisson}, it has for first consequence that $a_\ep$ is
globally compact in $L^p$, $(1-1/d)^{-1}<p<\infty$. 
Applying Th. \ref{mainth}, one deduces that $n_\ep$ is
locally compact in $L^1$ and by \eqref{queue} that $n_\ep$ is compact
in $L^1$ and so in any $L^p$, $1\leq p<\infty$.

Let us now extract two converging subsequences (still denoted by $\ep$) 
\[
a_\ep\longrightarrow a,\quad n_\ep\longrightarrow n.
\]
We may now easily pass to the limit in every term of \eqref{viscous}
and \eqref{poisson} to deduce that $n$ and $a$ are solutions, in the
sense of distributions, to \eqref{maineq} coupled with
\eqref{poisson}.

Proving that $n$ is an entropy solution to \eqref{maineq} follows the
usual procedure. For any $\phi\in C^2$ convex, we first note that
\[
\partial_t \phi(n_\ep)+{\rm div}_x\,(a_\ep q(n_\ep))+(\phi'(n_\ep) 
f(n_\ep)-q(n_\ep))\,{\rm div}_x\,a_\ep\leq 0,
\]
with $q'=\phi'\,f'$. With the compactness of $n_\ep$, one may pass to
the limit in each term and obtain the same property for $n$, which
concludes the proof of Corollary \ref{existence}.
%
\section{Proof of Theorem \ref{mainth}}
%
%
\subsection{The compactness criterion}
We first introduce the compactness criterion that we use. Define a family
$K_h(x)=1/(|x|+h)^d$ for $|x|\leq 1$ and $K_h$ non negative,
independent of $h$, with support in $B(0,2)$ and in
$C^\infty(\R^d\setminus B(0,1))$. 
\begin{lemma} A sequence of functions $u_k$, uniformly bounded in
  $L^p(\R^d)$ 
 is compact in $L^p_{loc}$ if 
\[
\limsup_k\; |\log h|^{-1}\; \int_{\R^{2d}}
K_h(x-y)\,|u_k(x)-u_k(y)|^p\,dx\,dy\longrightarrow
0\quad\mbox{as}\ h\rightarrow 0.
\]
Conversely if $u_k$ is globally compact in $L^p$ then the previous limit holds.
\label{compact}
\end{lemma}

\noindent{\bf Proof.} We recall that $u_k$ is  compact in $L^p$
iff 
\[
\delta(\eta)=
\eta^{-d}\sup_k \int_{|x-y|\leq \eta} |u_k(x)-u_k(y)|^p\,dx\,dy \longrightarrow
0\quad\mbox{as}\ \eta\rightarrow 0. 
\]
So assuming $u_k$ is compact, one simply decomposes
\[\begin{split}
&\sup_k\int_{\R^{2d}}
K_h(x-y)\,|u_k(x)-u_k(y)|^p\,dx\,dy\leq C\\
&\ +C\,\sum_{n\leq |\log h|}
\sup_k\int_{2^{-n-1}\geq|x-y|\leq 2^{-n}}
 2^{dn}|u_k(x)-u_k(y)|^p\,dx\,dy\\
&\quad\leq C+C\,\sum_{n\leq |\log h|}\delta(2^{-n}),
\end{split}\]
which gives the result.

Conversely  assume that
\[
\alpha(h)=\limsup_k |\log h|^{-1}\, \int_{\R^{2d}}
K_h(x-y)\,|u_k(x)-u_k(y)|^p\,dx\,dy\longrightarrow
0\ \mbox{as}\ h\rightarrow 0.
\]
Denote $\tilde K_h(x)=C_h |\log h|^{-1}\,K_h(x-y)$,
with $C_h$ s.t.
\[
\int \tilde K_h(x)\,dx=1,
\]
and therefore $\tilde K_h$ a convolution kernel. Note that $C_h$ is
bounded from below and from above uniformly in $h$. Now
\[\begin{split}
\|u_k-\tilde K_h\star_x u_k\|_{L^p}^p\leq & |\log h|^{-p}\int_{\R^d}\left(
\int_{\R^{d}} K_h(x-y) |u_k(x)-u_k(y)|\,dy\right)^p\,dy\\
\leq & |\log h|^{-p}\,\|K_h\|_{L^1}^{p-1}\,\int_{\R^{2d}} K_h(x-y)
|u_k(x)-u_k(y)|^p\,dy\,dx\\
\leq & C\,|\log h|^{-1}\, \int_{\R^{2d}} K_h(x-y)
|u_k(x)-u_k(y)|^p\,dy\,dx  
\end{split}\]
is converging to $0$ uniformly in $k$ as the $\limsup$ is $0$ and it is
converging for any fixed $k$ by the usual approximation by convolution
in $L^p$. On the other hand for a fixed
$h$, $\tilde K_h\star u_k$ is compact in $k$ and this proves that $u_k$ also is.
\subsection{The main argument given for a linear transport equation}
%
Before proving Theorem \ref{mainth}, we wish to explain the main idea
behind the proof in a  simple and wellknown setting. Let us consider a sequence
$u_\ep$ of solutions to the transport equation
\begin{equation}\begin{split}
&\partial_t u_\ep(t,x) +v_\ep\cdot\nabla u_\ep(t,x)=0,\quad t\in
        [0,\ T],\ x\in \R^d,\\
&u_\ep(t=0,x)=u_\ep^0(x),
\end{split}\label{transport}\end{equation}
for a given velocity field. The following result was originally proved
in \cite{DL}
\begin{theorem} Assume that $u_\ep^0$ is uniformly bounded in $L^1
\cap L^\infty$ and compact. 
Assume moreover that $v_\ep$ is compact in $L^p$,
uniformly bounded in $L^\infty_t
W^{1,p}_x$ for
some $p>1$ and that $\hbox{\em div}\,v_\ep=0$. Then the sequence of
solutions $u_\ep$ to \eqref{transport} is compact in $L^1$.
\label{comptransport}\end{theorem}
\noindent{\bf Proof of Theorem \ref{comptransport}.}

First of all notice that $u_\ep$ is uniformly bounded in
$L^\infty_t(L^1_x\cap L^\infty_x)$. Moreover as $v_\ep$ is compact,
one may freely assume that it converges toward a limit $v\in L^\infty_t
W^{1,p}_x$ (by extracting a subsequence).

Now define 
\[
Q_\ep(t)=\int_{\R^{2d}} K_h(x-y)\,|u_\ep(t,x)-u_\ep(t,y)|^2\,dx\,dy.
\] 
From Equation \eqref{transport} the divergence free condition on
$v_\ep$, 
one simply computes
\[
\frac{dQ_\ep}{dt}=\int_{\R^{2d}}
\nabla K_h(x-y)\,(v_\ep(t,y)-v_\ep(t,x))|u_\ep(t,x)-u_\ep(t,y)|^2\,dx\,dy. 
\]
Therefore by introducing the limit $v$ 
\[\begin{split}
\frac{dQ_\ep}{dt}\leq &C\,\|v_\ep-v\|_{L^p}\,\|\nabla K_h\|_{L^1}\\
&+\int_{\R^{2d}}
\nabla K_h(x-y)\,(v(t,y)-v(t,x))|u_\ep(t,x)-u_\ep(t,y)|^2\,dx\,dy.
\end{split}\]
The second term is equal to 
\[
\int_0^1\int_{\R^{2d}}
(x-y)\otimes\nabla K_h(x-y):\nabla v(t,\theta x+(1-\theta)y)\,
|u_\ep(t,x)-u_\ep(t,y)|^2\,dx\,dy,
\]
with $A:B$ denoting the full contraction of the two matrices. Note
that for $|x|>1$, $\nabla K_h$ is bounded and for $|x|<1$,
\[
x\otimes \nabla K_h(x)=\frac{x\otimes x}{(|x|+h)^{d+1}|x|}.
\]
Define
\[
\tilde K_h(x)=x\otimes \nabla K_h(x)-\lambda
\frac{\hbox{Id}\,|x|}{(|x|+h)^{d+1}}\,{\mathbb I}_{|x|\leq 1},
\]
with $\lambda=\int_{S^{d-1}} \omega_1^2\,d\omega$.

Thanks to the definition of $\lambda$, $\tilde K_h$ is now a Calderon-Zygmund
operator, meaning that for any $1<q<\infty$, there exists a constant
$C$ independent of $h$ s.t.
\[
\|\tilde K_h\star g\|_{L^q}\leq C\,\|g\|_{L^q}.
\]
As $v$ is divergence free, one may simply replace by $\tilde K_h$ 
\[\begin{split}
&\int_0^1\int_{\R^{2d}}
(x\!-\!y)\otimes\nabla K_h(x\!-\!y):\nabla v(t,\theta x+(1\!-\!\theta)y)\,
|u_\ep(t,x)-u_\ep(t,y)|^2\,dx\,dy\\
&\quad=\int_0^1\int_{\R^{2d}}
\tilde K_h(x-y):\nabla v(t,\theta x+(1-\theta)y)\,
|u_\ep(t,x)-u_\ep(t,y)|^2\,dx\,dy\\
&\quad\leq C\,\int_0^1\int_{\R^{2d}}
K_h(x-y)\,|\nabla v(t,\theta x+(1-\theta)y)-\nabla v(t,x)|\,dx\,dy\\
&\qquad\ + \|\tilde K_h\star (\nabla v u_\ep^2)\|_{L^1}+
2\,\|u_\ep\,\tilde K_h\star (\nabla v u_\ep)\|_{L^1}+\|u_\ep^2\,\tilde
K_h\star
\nabla v\|_{L^1}.
\end{split}\]
Thanks to the uniform bounds on $u_\ep$, and changing variables, one
immediately deduce that 
\[\begin{split}
&\int_0^1\int_{\R^{2d}}
(x\!-\!y)\otimes\nabla K_h(x\!-\!y):\nabla v(t,\theta x+(1\!-\!\theta)y)\,
|u_\ep(t,x)-u_\ep(t,y)|^2\,dx\,dy\\
&\ \leq C+C\,\int_{\R^{2d}}
K_h(x-y)\,|\nabla v(t,x)-\nabla v(t,y)|\,dx\,dy.
\end{split}\]
Putting together all the terms in the estimate, we have
\[\begin{split}
&\frac{dQ_\ep}{dt}\leq C+C\,\frac{\|v_\ep-v\|_{L^p}}{h}\\
&\quad +C\,\int_{\R^{2d}}
K_h(x-y)\,|\nabla v(t,x)-\nabla v(t,y)|\,dx\,dy
\end{split}\]
or
\[\begin{split}
Q_\ep(t)\leq &C+C\,\frac{\|v_\ep-v\|_{L^p}}{h}\\
&+C\,\int_0^T\int_{\R^{2d}}
K_h(x-y)\,|\nabla v(t,x)-\nabla v(t,y)|\,dx\,dy\\
&+C\,\int_{\R^{2d}}
K_h(x-y)\,|n_\ep^0(x)-n_\ep^0(y)|^2\,dx\,dy.
\end{split}\]
As $n_\ep^0$ is compact and $v$ is independent of $\ep$ then the
previous estimate shows that
\[
\lim_{h\rightarrow 0}|\log h|^{-1}\limsup_\ep \sup_t\int_{\R^{2d}}
K_h(x-y)\,|n_\ep(t,x)-n_\ep(t,y)|^2\,dx\,dy=0.
\] 
Lemma \ref{compact} then proves that $u_\ep$ is compact in
space. However by Eq. \eqref{transport}, $\partial_t u_\ep$ is
uniformly bounded in $L^\infty_t(W^{-1,p}_x)$. Therefore compactness
in time follows and the theorem is proved.
%
\subsection{A simple proof for Theorem \ref{mainth}}
%
We first give here a simple proof of the compactness. This proof is
not optimal in the sense that it does not give an explicit rate for
how the norm in our compactness criterion behaves
\[
\int_0^T\int K_h(x-y) |n_\ep(t,x)-n_\ep(t,y)|\,dx\,dy\,dt.
\]
This is however a more difficult problem, which is partially dealt
with in the next section.

As $a_\ep$ is
compact in $L^p$, by extracting a subsequence (still denoted by $\ep$),
$a_\ep$ 
converges strongly in $L^p$ to some $a\in W^{1,p}$. By the compactness
of $d_\ep$ and $n^0_\ep$ and by Lemma \ref{compact},
we may assume without loss
of generality that there exists a continuous function $\delta(h)$ with
$\delta(0)=0$, independent of $\ep$ and a function $\alpha(\ep)$, 
 s.t. 
\begin{equation}\begin{split}
& |\log h|^{-1}\;
\int_{\R^d} K_h(x-y)\,|n^0_\ep(y)-n^0_\ep(x)|\,dx\,dy\leq \delta(h),\\
&
|\log h|^{-1}\;\int_0^T\int_{\R^d} K_h(x-y)\,
|d_\ep(t,y)-d_\ep(t,x)|
\,dx\,dy\,dt\leq \delta(h), \\
&
|\log h|^{-1}\;\int_0^T\int_{\R^d} K_h(x-y)\,
|\nabla a(t,y)-\nabla a(t,x)|^p
\,dx\,dy\,dt\leq \delta^p(h), \\
&\int_0^T\int_{\R^d} 
|a_\ep(t,x)-a(t,x)|^p\;dx\,dt\leq \alpha^p(\ep).
\end{split}
\label{convergea}
\end{equation}
Note that the estimate is written for $\nabla a$ and not for the
sequence $\nabla a_\ep$ as no compactness can be assumed on $\nabla a_\ep$.

Then one proves
\begin{prop} 
Let $n_\ep$ be a sequence of solutions to \eqref{viscous} with initial data
$n_\ep^0$ uniformly bounded in $L^1\cap L^\infty$ and compact in
$L^1$. Assume
\eqref{w11}, \eqref{divb}, \eqref{divcomp} and hence
\eqref{convergea}. Then for some constant $C$ uniform in $h$ and $\ep$
\[
\int_0^T \int_{\R^{2d}} K_h(x-y)\,|n_\ep(t,x)-n_\ep(t,y)|\,
dx\,dy\,dt\leq C\,\frac{\ep^2}{h^2}+
C\,\delta(h)\,|\log h|+C\,\frac{\alpha(\ep)}{h}.
\]\label{propsimple}
\end{prop}
The disappointing part of Prop. \ref{propsimple} is that the rates
$\delta(h)$ and $\alpha(\ep)$ are not explicit but depend
intrinsically on the sequence $a_\ep$. See the next section for a more
explicit (but much more complicated) result.

Prop. \ref{propsimple} proves the compactness in space of $n_\ep$  by
Lemma \ref{compact}. The compactness in time is then straightforward
since $n_\ep$ solves a transport equation \eqref{maineq}.

Hence Theorem \ref{mainth} follows.

{\noindent \bf Proof of Prop. \ref{propsimple}.}

The proof mostly follows the steps of the proof of Theorem
\ref{comptransport}. The main differences are the nonlinear flux, the
vanishing viscosity terms and
the fact that now the field
$a_\ep$ is not assumed to be divergence free (only bounded).

First of all, by 
condition \eqref{divb}, for any $T>0$, $n_\ep(t,x)$
is bounded in $L^1\cap L^\infty([0,\ T]\times \R^d)$, uniformly in $\ep$.

We start with Kruzkov's usual argument of doubling of variable. If
$n_\ep$ is a solution to \eqref{viscous} then
\[\begin{split}
&\partial_t |n_\ep(t,x)-n_\ep(t,y)|+\hbox{div}_x
\big(a_\ep(t,x)F(n_\ep(t,x),n_\ep(t,y))\big)\\
&\quad +\hbox{div}_y
\big(a_\ep(t,y)F(n_\ep(t,y),n_\ep(t,x))\big)+\hbox{div}_x a_\ep(t,x)\;
G(n_\ep(t,x),n_\ep(t,y)) \\
&\quad+\hbox{div}_y a_\ep(t,y)\;
G(n_\ep(t,y),n_\ep(t,x))-\ep^2\,(\Delta_x+\Delta_y)\,|n_\ep(t,x)-n_\ep(t,y)|
\leq 0. 
\end{split}\] 
This computation is formal but can easily be made rigourous by using a
suitable regularisation of $|.|$. Here $F$ satisfies
\[
F'(\xi,\zeta)=f'(\xi)\,\hbox{sign}(\xi-\zeta),\ F(\xi,\zeta)=0,
\]
which means that
\[
F(\xi,\zeta)=(f(\xi)-f(\zeta))\,\hbox{sign}(\xi-\zeta)=F(\zeta,\xi).
\]
And as for $G$
\[
G(\xi,\zeta)=f(\xi)\,\hbox{sign}(\xi-\zeta)-F(\xi,\zeta)=-G(\zeta,\xi).
\]
Now define 
\[
Q(t)=\int_{\R^{2d}} K_h(x-y)\,|n_\ep(t,x)-n_\ep(t,y)|\,dx\,dy.
\]
Remark that
\[\begin{split}
\ep^2\int_{\R^{2d}} &K_h(x-y)\,(\Delta_x+\Delta_y)\,
|n_\ep(t,x)-n_\ep(t,y)|\,dx\,dy\\
&=\ep^2\int_{\R^{2d}}  \Delta K_h(x-y)\,
|n_\ep(t,x)-n_\ep(t,y)|\,dx\,dy\\
&\leq C\,\ep^2 \|\Delta K_h\|_{L^1}\leq C\,\frac{\ep^2}{h^2}.
\end{split}\]
Using this and because of the symmetry of $F$ and the
antisymmetry of $G$
\[\begin{split}
\frac{d}{dt}Q(t) 
\leq & C\,\frac{\ep^2}{h^2}+
\int_{|x-y|\leq 1} \frac{x-y}{(|x-y|+h)^{d+2}}\cdot
(a_\ep(t,y)-a_\ep(t,x))\\
&\qquad\qquad\qquad\qquad\qquad\qquad\qquad
F(n_\ep(t,y),n_\ep(t,x))\,dx\,dy\\
+&\int_{|x-y|\geq 1} \nabla K(x-y)\cdot (a_\ep(t,y)-a_\ep(t,x))\,
F(n_\ep(t,y),n_\ep(t,x))\,dx\,dy\\
+&\int_{\R^{2d}} K_h(x-y)\,(\hbox{div}\,
a_\ep(t,x)-\hbox{div}\,
a_\ep(t,y))\,G(n_\ep(t,y),n_\ep(t,x))\,dx\,dy\\
=&C\,\frac{\ep^2}{h^2}+A+B+D. 
\end{split}\]
Let us begin with the last term. Use \eqref{divcomp} to decompose
\[\begin{split}
D\leq&\int_{\R^{2d}} K_h(x-y)\,|d_\ep(t,x)-d_\ep(t,y)|\;
|G(n_\ep(t,y),n_\ep(t,x))|\,dx\,dy\\
&+C \int_{\R^{2d}} K_h(x-y)\,|n_\ep(t,x)-n_\ep(t,y)|\;
|G(n_\ep(t,y),n_\ep(t,x))|\,dx\,dy.
\end{split}
\]
As $G(n_\ep(t,x),n_\ep(t,y))$ is uniformly bounded in $L^\infty$, one
gets from \eqref{convergea}
\begin{equation}
\int_0^T D\,dt\leq  |\log h|\,\delta(h)+C\,\int_0^T Q(t)\,dt.
\end{equation}
For the second term $B$, just note that $\nabla K\in
C^\infty_c(\R^d\setminus B(0,1))$, and that
$|F(n_\ep(t,x),n_\ep(t,x))|\leq |f(n_\ep(t,x))|+|f(n_\ep(t,y))|$ is
uniformly bounded in $L^1$. So one simply has
\[
B\leq C.
\]
The main term is hence $A$. Using again \eqref{convergea} and the
bound on $|F(.,.)|$, one gets
\[\begin{split}
\int_0^T A\,dt
&\leq C\,\frac{\alpha(\ep)}{h}
+\int_0^T\int_{|x-y|\leq 1} \frac{x-y}{(|x-y|+h)^{d+1}\,|x-y|}\cdot
(a(t,y)-a(t,x))\\
&\qquad\qquad\qquad\qquad\qquad\quad 
F(n_\ep(t,y),n_\ep(t,x))\,dx\,dy\,dt\\
\leq C\,\frac{\alpha(\ep)}{h}&+
\int_0^T\int_0^1\int_{|x-y|\leq 1} 
\frac{(x-y)\otimes(x-y)}{(|x-y|+h)^{d+1}\,|x-y|}\,
:\nabla a(t,\theta x+(1-\theta)y)\\
&\qquad\qquad\qquad\qquad\qquad\qquad\qquad 
F(n_\ep(t,y),n_\ep(t,x))\,dx\,dy\,d\theta\,dt.
\end{split}
\]
Still using \eqref{convergea},
\[
\begin{split}
&\int_0^T A\,dt
\leq C\,\left(\frac{\alpha(\ep)}{h}+|\log h|\,\delta(h)\right)\\
&+
\int_0^T\int_{|x-y|\leq 1} \frac{(x-y)\otimes(x-y)}{(|x-y|+h)^{d+1}\,|x-y|}\,
:\nabla a(t,x)
F(n_\ep(t,y),n_\ep(t,x))\,dx\,dy\,dt\\
&\qquad\qquad
\leq C\,\left(\frac{\alpha(\ep)}{h}+|\log h|\,\delta(h)\right)
+\int_0^T E(t)\,dt.
\end{split}
\]
Denote as in the proof of Theorem \ref{comptransport}
\[
\lambda=\int_{S^{d-1}} \omega_1^2\,dS(\omega),\quad \bar
K_h(x)=\left(\frac{x\otimes x}{(|x|+h)^{d+1}\,|x|}-\lambda
\frac{|x|}{(|x|+h)^{d+1}}\,Id\right)\;{\mathbb I}_{|x|\leq 1}. 
\]
By the definition of $\lambda$, $\bar K_h$ is a Calderon-Zygmund
operator bounded on any $L^p$ for $1<p<\infty$.
Now write
\[\begin{split}
E&=\int_{\R^{2d}} \bar K_h(x-y)\,\nabla
a(x)\,F(n_\ep(t,y),n_\ep(t,x))\,dx\,dy \\
&\quad
+\lambda\,\int_{|x-y|\leq
  1}\frac{|x-y|}{(|x-y|+h)^{d+1}}\,\hbox{div}\,a(t,x)
\;F(n_\ep(t,y),n_\ep(t,x))\,dx\,dy\\
&\leq \int_{\R^{2d}} \bar K_h(x-y)\,\nabla
a(x)\,F(n_\ep(t,y),n_\ep(t,x))\,dx\,dy+C\,Q(t),
\end{split} \]
as the divergence of $a$ is bounded.

Introduce 
\[
\chi_\ep(t,x,\xi)={\mathbb I}_{0\leq \xi\leq n_\ep(t,x)}.
\]
Then note that $\chi_\ep$ is compactly supported in $\xi$ and that
\[
F(n_\ep(t,y),n_\ep(t,x))=\int_0^\infty
f'(\xi)\,|\chi_\ep(t,x,\xi)-\chi_\ep(t,y,\xi)|^2\,d\xi.
\]
Hence as $\nabla a\in L^p$, and $\chi_\ep$ is uniformly bounded in
$L^\infty_{t,\xi}(L^1_x\cap L^\infty_x)$, for $1/p+1/p^*=1$,
\[\begin{split}
\int_{\R^{2d}} &\bar K_h(x-y)\,\nabla
a(x)\,F(n_\ep(t,y),n_\ep(t,x))\,dx\,dy\\
&=\int_{\R_+} f'(\xi)\,\int_{\R^{2d}} \bar
K_h(x-y)\,\nabla 
a(x)\,|\chi_\ep(t,x,\xi)-\chi_\ep(t,y,\xi)|^2\,d\xi\,dx\,dy\\
&\leq \int_{\R_+} f'(\xi)\,\big(\|\bar K_h\star (\nabla a
\chi_\ep^2)\|_{L^1}+\|\bar K_h\star 
\chi_\ep^2)\|_{L^{p^*}}+2\|\bar K_h\star (\nabla a
\chi_\ep)\|_{L^1}\big)\,d\xi\\
&\leq C.
\end{split}\]

Combining all estimates we conclude that
\[
Q(t)\leq Q(0)+C\,\frac{\ep^2}{h^2}
+C\,|\log h|\,\delta(h)+C\,\frac{\alpha(\ep)}{h}+C\,\int_0^t Q(s)\,
ds. 
\]
The initial data $Q(0)$ is bounded by \eqref{convergea} and finally by
Gronwall lemma we obtain on any finite interval
\[
Q(t)\leq C\,\frac{\ep^2}{h^2}+C\,|\log h|\,\delta(h)+C\,\frac{\alpha(\ep)}{h},
\]
which proves the proposition.
%
\section{An explicit estimate : Proof of Theorem \ref{quantth}}
%
Checking carefully the proof of Prop. \ref{propsimple}, one sees that 
to get an explicit rate, it would be necessary to bound a term like
\begin{equation}
\int_{\R^{2d}} \nabla
K_h(x-y)\,(a_\ep(x)-a_\ep(y))\,|g_\ep(x)-g_\ep
(y)|^2\,dx\,dy\label{calderon}
\end{equation}
only in terms of the $W^{1,p}$ norm of $a_\ep$ and the $L^1\cap
L^\infty$ norms of $g_\ep$. 

Here we do not aim at optimal estimates, just explicit ones. 
We present a very elementary
proof of 
\begin{prop} Let $1<p<\infty$,  $\exists C_p<\infty$ s.t. $\forall
  a(x),\ g(x)$ smooth and compactly supported
\[\begin{split}
\int_{\R^{2d}} &\nabla K_h(x-y)\,(a(x)-a(y))\,|g(x)-g(y)|^2\,dx\,dy \\
&\leq C_p\,\|g\|_{L^\infty}\,\|g\|_{L^1\cap L^{p^*}}\,\|\nabla
a\|_{L^p\cap L^1}
\,|\log
h|^{1/\bar p}\\
&+C_p\,(\|\hbox{\em div}\, a\|_{L^\infty}+\|\nabla a\|_{L^p})\,\int_{\R^{2d}} 
K_h(x-y)|g(x)-g(y)|^2\,dx\,dy,
\end{split}\]
with $1/p^*+1/p=1$ and $\bar p=\min(p,2)$.
\label{theprop}
\end{prop}
Note that the rate $|\log h|^{1/\bar p}$ is most probably not
optimal. A way to obtain a better rate could be to combine
Lemma \ref{teclemma} below with the estimates in \cite{MTT2} as we
suggest below.

The kind of Calderon-like estimate like Prop. \ref{theprop} 
has been extensively studied in
dimension $1$, see for instance \cite{Ca} or \cite{CM},
\cite{CM2}. The situation in higher dimension is however more
complicated. In particular it seems necessary to use the bound on the
divergence of $a_\ep$ to estimate \eqref{calderon} (as was already
suggested by the proof of Prop. \ref{propsimple}).

Following the previous section, a simple idea would be to estimate
\eqref{calderon} by
\[\begin{split}
& C\,\|\hbox{div}\, a\|_{L^\infty}\,\int_{\R^{2d}}
K_h(x-y)\,|g_\ep(x)-g_\ep(y)|^2\,dx\,dy\\
&+\int_0^1\int_{\R^{2d}}
L_h(x-y):\nabla a(\theta x+(1-\theta)y)\,|g_\ep(x)-g_\ep(y)|^2\,dx\,dy,
\end{split}\]  
where $L_h$ is now a Calderon-Zygmund operator. Expanding the square,
one sees that it would be enough to bound in some $L^q$ space
\[
\int_0^1\int_{\R^{d}}
L_h(x-y):\nabla a(\theta x+(1-\theta)y)\,g_\ep(y)\,dy.
\] 
Using Fourier transform (we denote by ${\cal F}$ the Fourier
transform)
and an easy change of variable, this term is
equal to
\[
\int_{\R^{2d}} e^{ix\cdot (\xi_1+\xi_2)}\, m(\xi_1,\xi_2)\,{\cal
    F}\nabla a (\xi_1)\,{\cal F} g(\xi_2)\,d\xi_1\,d\xi_2,
\]
with
\[
m(\xi_1,\xi_2)=\int_0^1 {\cal F}L_h(\theta\,\xi_1+\xi_2)\,d\theta.
\]
We now have a multi-linear operator in dimension $d$ of the kind
studied in Muscalu, Tao, Thiele 
\cite{MTT2}. Unfortunately $m$ does not satisfy the
assumptions of this last article as it does not have the right behaviour on
the subspace $\xi_1\parallel \xi_2$. Instead it would be necessary to
have a multi-dimensional equivalent of \cite{MTT} (which, as far as we
know, is not yet proved) or to use Lemma \ref{teclemma}. 
%
\subsection{Proof of Theorem \ref{quantth} given Prop. \ref{theprop}}
%
For the moment let us assume Prop. \ref{theprop}. Define
\[
Q(t)=\int_{\R^{2d}} K_h(x-y)\,|n_\ep(t,x)-n_\ep(t,y)|\,dx\,dy.
\]
We follow the same first steps as in the proof of
Prop. \ref{propsimple}, with the same notations. We obtain
\[\begin{split}
\frac{dQ}{dt}\leq &C+C\,\frac{\ep^2}{h^2}+C\,Q(t)+C\,\int_{\R^{2d}}
K_h(x-y)\,|d_\ep(t,x)-d_\ep(t,y)|\,dx\,dy \\
&+\int_{\R^{2d}} \nabla K_h(x-y)\cdot
(a_\ep(t,x)-a_\ep(t,y))\,F(n_\ep(t,y),n_\ep(t,x))\,dx\,dy.
\end{split}\]
We only have to bound the last term. Let us introduce again
\[
\chi_\ep(t,x,\xi)={\mathbb I}_{0\leq \xi\leq n_\ep(t,x)}.
\]
Note that $\chi_\ep$ is supported in $\xi$ in
$[0,\ \|n_\ep^0\|_{L^\infty}]\subset [0,\ C]$. 

Now write
\[\begin{split}
&\int_{\R^{2d}} \nabla K_h(x-y)\cdot
(a_\ep(t,x)-a_\ep(t,y))\,F(n_\ep(t,y),n_\ep(t,x))\,dx\,dy\\
&\quad=\int_0^C f'(\xi)
\int_{\R^{2d}} \nabla K_h(x-y)\cdot
(a_\ep(t,x)-a_\ep(t,y))\\
&\qquad\qquad\qquad\qquad\qquad\qquad
|\chi_\ep(t,y,\xi)-\chi_\ep(t,x,\xi)|^2\,dx\,dy\,d\xi\\
&\quad\leq C\, |\log h|^{1-2/\bar p}+C\,\int_{\R^{2d}} 
K_h(x-y)\,\int_0^C |\chi_\ep(t,y,\xi)-\chi_\ep(t,x,\xi)|^2\,d\xi\,dx\,dy,  
\end{split}\]
using Prop. \ref{theprop} and the uniform bounds on
$\|a_\ep\|_{L^\infty_tL^p_x}$ and $\|\chi_\ep\|_{L^1\cap L^\infty}$. 
Now simply note that because of the
definition of $\chi_\ep$
\[
\int_0^C |\chi_\ep(t,y,\xi)-\chi_\ep(t,x,\xi)|^2\,d\xi\leq
|n_\ep(t,x)-n_\ep(t,y)|, 
\]
and the last term in the previous inequality is hence simply bounded
by $Q$.
One finally obtains
\[\begin{split}
\frac{dQ}{dt}\leq &C++C\,\frac{\ep^2}{h^2}+
C\,Q(t)+C\, |\log h|^{1-2/\bar p}\\
&+C\,\int_{\R^{2d}}
K_h(x-y)\,|d_\ep(t,x)-d_\ep(t,y)|\,dx\,dy. 
\end{split}\]
To conclude the proof of Theorem \ref{quantth}, it is now enough to
apply Gronwall's lemma.
%
\subsection{ Beginning of the proof of Prop. \ref{theprop}\label{begproof}} 
%
As before we will control
$a(x)-a(y)$ with $\nabla a$. Contrary to the previous case though, it
is not enough to integrate over the segment. Instead use the lemma
\begin{lemma}
\[\begin{split}
a_i(x)-a_i(y)&=|x-y|\,\int_{B(0,1)}
\psi\left(z,\frac{x-y}{|x-y|}\right)\cdot \nabla
a_i\left(x+|x-y|\,z\right)\frac{dz}{|z|^{d-1}}\\
&+|x-y|\,\int_{B(0,1)}
\psi\left(z,\frac{x-y}{|x-y|}\right)\cdot \nabla
a_i\left(y+|x-y|\,z\right)\frac{dz}{|z|^{d-1}},
\end{split}\]
where $|z|\,\psi$ is Lipschitz on $B(0,1)\times S^{d-1}$  
and for a  given
constant $\alpha$, 
\[
\int_{B(0,1)}
\psi\left(z,\frac{x-y}{|x-y|}\right)\,\frac{dz}{|z|^{d-1}}
=\alpha
\frac{x-y}{|x-y|}.
\]\label{teclemma}
\end{lemma}
\noindent{Proof of Lemma \ref{teclemma}.} We refer to \cite{CJ} for a
complete, detailed proof. Let us simply mention that the idea is to
integrate along many trajectories between $x$ and $y$ instead of just
the segment.\cqfd

\bigskip

Lemma \ref{teclemma} gives two terms that are completely symmetric and
it is enough to deal with one of them. After an easy change of
variable, one finds
\[\begin{split}
\int_{\R^{2d}} &\nabla K_h(x-y)\,(a(x)-a(y))\,|g(x)-g(y)|^2\,dx\,dy\\
&=\int_{\R^{d}}\int_0^1 \frac{r^{d}}{(r+h)^{d+1}}
\int_{B(0,1)}\int_{S^{d-1}} \psi(z,\omega)\otimes \omega\,:\,
\nabla a(x+rz)\\
&\qquad\qquad |g(x)-g(x+r\omega)|^2\,d\omega\,\frac{dz}{|z|^{d-1}}
\,dr\,dx+\quad symmetric, 
\end{split}\]
where $A:B$ denotes the total contraction of two matrices $\sum_{i,j}
A_{ij}\, B_{ij}$.

Now define
\[
L(z,\omega)=\psi(z,\omega)\otimes
\omega-\lambda\,\hbox{Id}, 
\]
 for
 $\lambda=\int_{B(0,1)}\int_{S^{d-1}}\omega_1^2\,d\omega\,
\frac{dz}{|z|^{d-1}}$. 

Note that
\[\begin{split}
&\int_{\R^{d}}\int_0^1 \frac{r^{d}}{(r+h)^{d+1}}
\int_{B(0,1)}\int_{S^{d-1}} \psi(z,\omega)\otimes \omega\,:\,
\nabla a(x+rz)\\
&\qquad\qquad |g(x)-g(x+r\omega)|^2\,d\omega\,\frac{dz}{|z|^{d-1}}\,dr\,dx\\
&\quad\leq \int_0^1\int_{\R^d\times B(0,1)\times S^{d-1}}
\frac{r^{d}\,
L(z,\omega)}{(r+h)^{d+1}}\,:\,
\nabla a(x+rz)\,|g(x)-g(x+r\omega)|^2\,\frac{1}{|z|^{d-1}}\\
&\qquad+C\,\|\hbox{div}\, a\|_{L^\infty} \int_0^1\int_{\R^d\times
  S^{d-1}}\frac{r^{d}}{(r+h)^{d+1}}\, |g(x)-g(x+r\omega)|^2
\end{split}\]
 By the definition of $K_h$, the second term is
bounded by
\[
C\,\|\hbox{div}\, a\|_{L^\infty}\int_{\R^{2d}} 
K_h(x-y)|g(x)-g(y)|^2\,dx\,dy.
\]
and it only remains to bound the first one. In order to get the
optimal rate for $\nabla a\in L^p$ with $p>2$, we need to introduce an
additional decomposition of $\nabla a$. For $p>2$ as $L^p$ may be
obtained by interpolating between $L^2$ and $L^\infty$, let
\[
\nabla a=A+\bar A,\quad \|\bar A\|_{L^\infty}\leq 2\|\nabla
a\|_{L^p},\quad
\|A\|_{L^2}\leq 2 \|\nabla a\|_{L^p}.
\] 
If $p<2$ then we simply put $A=\nabla a$. In both cases, if $\nabla a$
is smooth and compactly supported then one may of course assume the
same of $A$ and $\bar A$.

Define
\[
Q(A,g)=\int_0^1\int_{B(0,1)\times S^{d-1}}\frac{r^{d}\,
L(z,\omega)}{(r+h)^{d+1}}\,:\,
A(x+rz)\,g(x+r\omega)\,d\omega\,\frac{dz}{|z|^{d-1}}\,dr. 
\]
The term with $\bar A$ may be bounded directly by using the $L^\infty$
norm of $\bar A$; for the other one 
simply by expanding the square $|g(x)-g(y)|^2$, one obtains
\[
\begin{split}
\int_{\R^{2d}} &\nabla K_h(x-y)\,(a(x)-a(y))\,|g(x)-g(y)|^2\,dx\,dy\\
&\leq C\,(\|\hbox{div}\, a\|_{L^\infty}+\|\nabla a\|_{L^p})
\, \int_{\R^{2d}}K_h(x-y)\, 
|g(x)-g(y)|^2\,dx\,dy\\
&+\int_{\R^d} (-2g\,Q(A,g)+g^2\,Q(A,1))\,dx.
\end{split}
\]
Note that bounding $Q(A,1)$ is in fact easy as it is an ordinary
convolution and
 $\frac{1}{r}\, L$ defines a Calderon-Zygmund operator. However the control of
$Q(A,g)$ essentially requires to rework Calderon-Zygmund theory.

Of course for $r$ of order $h$ then one has
\[\begin{split}
&\left\|\int_0^h\int_{B(0,1)\times S^{d-1}}\frac{r^{d}\,
L(z,\omega)}{(r+h)^{d+1}}\,:\,
A(x+rz)\,g(x+r\omega)\,d\omega\,\frac{dz}{|z|^{d-1}}\,dr\right\|_{L^1}\\
&\qquad \leq \frac{1}{h}\int_0^h\int_{B(0,1)\times S^{d-1}}
\left\|A(x+rz)\,g(x+r\omega)\right\|_{L^1}\,d\omega\,\frac{dz}{|z|^{d-1}}\,dr\\
&\qquad\leq \|A\|_{L^1}\,\|g\|_{L^\infty}.\end{split}\]

It is hence enough to consider
\[
\bar Q(A,g)=\int_h^1\int_{B(0,1)\times S^{d-1}}\frac{r^{d}\,
L(z,\omega)}{(r+h)^{d+1}}\,:\,
A(x+rz)\,g(x+r\omega)\,d\omega\,\frac{dz}{|z|^{d-1}}\,dr. 
\]
We introduce the Littlewood-Paley decomposition of $A$
(see for instance Triebel \cite{Tr})
\[
A(x)= \sum_{i=0}^\infty A_i(x),
\] 
where for $i>0$, $\hat A_i={\cal F}{A}\, p(2^{-i}\xi)$ with $p$
compactly supported in the annulus of radii $1/2,\;2$; and $\hat
A_0={\cal F}{A}\, p_0(\xi)$. The functions $p$ and $p_0$
determines a partition of unity. We note either ${\cal F}g$ or $\hat
g$ the Fourier transform of a function $g$. In the following we denote
by ${\cal P}_i$ the projection operator ${\cal P}_i \phi={\cal F}^{-1}
p(2^{-i}\xi)\,\hat \phi$.

There is an obvious critical scale in the decomposition which is
where $2^{-i}$ is of order $r$. Accordingly we decompose further
\[\begin{split}
&\bar Q(A,g)= Q_1(A,g)+ Q_2(A,g)\\
&\ =\sum_{i\leq |\log h|} \int_h^{2^{-i}} \int_{B(0,1)\times S^{d-1}}\frac{r^{d}\, 
L(z,\omega)}{(r+h)^{d+1}}\,:\,
A_i(x+rz)\,g(x+r\omega)\,d\omega\,\frac{dz}{|z|^{q-1}}\,dr\\
&\ +\sum_{i} \int_{\max(h,2^{-i})}^1 \int_{B(0,1)\times S^{d-1}}\frac{r^{d}\, 
L(z,\omega)}{(r+h)^{d+1}}\,:\,
A_i(x+rz)\,g(x+r\omega)\,d\omega\,\frac{dz}{|z|^{q-1}}\,dr.
\end{split}\] 
Each term is bounded in a different way. Note of course that in $Q_1$
as $r\geq h$ there is of course no frequency $i$ higher than $|\log
h|$ (they are all in $Q_2$). 
%
\subsection{Control on $Q_1$ in $L^2$}
The aim is here is to prove
\begin{lemma}  
$\forall 1<q<\infty$, $\exists C>0$ such that for any $A$ and $g$
  smooth and compactly supported 
  functions, 
\[
\|Q_1(A,g)\|_{L^1}\leq C\,\|A\|_{B^0_{q,1}}\,\|g\|_{L^{q^*}}.
\]
where $B^0_{q,1}$ is the usual Besov space and $1/q^*+1/q=1$. 
\label{Q1}
\end{lemma}
As we wish to remain as elementary as possible here, we avoid the use
of Besov spaces in the sequel. Instead for $q=2$ it it possible to
obtain directly the Lebesgue space by losing $|\log h|^{1/2}$ namely
\begin{lemma}  
 $\exists C>0$ such that for any $A$ and $g$
  smooth and compactly supported 
  functions, 
\[
\|Q_1(A,g)\|_{L^1}\leq C\,|\log h|^{1/2}\,\|g\|_{L^2}\,\|A\|_{L^{2}}.
\]\label{Q1l2}
\end{lemma} 

The proof is relatively simple. Indeed in $Q_1$ since $r<2^{-i}$,
$A_i$ does not change much over a ball of radius $r$. Therefore, we
simply replace $A_i(x+rz)$ by $A_i(x)$ in $Q_1$. This gives
\begin{equation}\begin{split}
&Q_1(A,g)\leq I+II\\
&\ \leq\sum_{i\leq |\log h|} \int_h^{2^{-i}} \int_{B(0,1)\times S^{d-1}}\frac{r^{d}\, 
L(z,\omega)}{(r+h)^{d+1}}\,:\,
A_i(x)\,g(x+r\omega)\,d\omega\,\frac{dz}{|z|^{d-1}}\,dr\\
&\quad+\sum_{i\leq |\log h|}\, \int_h^{2^{-i}}\frac{1}{r+h}\int_{S^{d-1}}
|g(x+r\omega)|\,d\omega\\
&\qquad\qquad\qquad\qquad\int_{B(0,1)}
|A_i(x+rz)-A_i(x)|\,\frac{dz}{|z|^{d-1}}\,dr.
\end{split}\label{decompQ1}
\end{equation}
Let us bound the first term. As $A_i$ does not depend on $z$ anymore,
this term is simply equal to
\[
\sum_{i,j\leq i} A_i(x) \int_h^{1} 
\int_{S^{d-1}}\tilde L_i(r\omega)
g_j(x+r\omega)\,r^{d-1}d\omega\,dr,
\]
where 
\[
\tilde
L_i(r\omega)=\frac{r}
{(r+h)^{d+1}}\,(\omega\otimes\omega-\tilde\lambda I)=\int_{B(0,1)}
\frac{r{{\mathbb I}_{r\leq 2^{-i}}}}
{(r+h)^{d+1}}\;L(z,\omega)\,\frac{dz}{|z|^{d-1}}. 
\]
By the definition of $\lambda$, $\tilde \lambda=\int_{S^{d-1}}
\omega_1^2\,d\omega$ and hence $\tilde 
L_i$ is a Calderon-Zygmund operator with operator norm bounded 
uniformly in $i$.

Now write for $1/q^*+1/q=1$
\[\begin{split}
\|I(x)\|_{L^1}
&\leq \sum_{i\leq |\log h|} 
\|A_i\|_{L^q}\,\|\tilde L_i\star g\|_{L^{q^*}}\\
&\leq C\,\|g\|_{L^{q^*}}\,\sum_i
\|A_i\|_{L^q}=C\,\|g\|_{L^{q^*}}\,\|A\|_{B^0_{q,1}}. 
\end{split}
\]

\medskip

Let us turn to the second term. 
\[\begin{split}
\|II\|_{L^1}&\leq \sum_{i\leq |\log h|}
\int_{S^{d-1}}\int_h^{2^{-i}}\left\|
 \frac{g(x+r\omega)}{r+h}
\int_{B(0,1)}
|A_i(x+rz)-A_i(x)|\,\frac{dz}{|z|^{d-1}}\right\|_{L^1}\\
&\leq C\,\|g\|_{L^{q^*}}\,\sum_{i\leq |\log h|} \int_{B(0,1)}\int_{h}^{2^{-i}}
 \|A_i(.+rz)-A_i(.)\|_{L^q}\;\frac{dr}{r+h}\,\frac{dz}{|z|^{d-1}}.\\
\end{split}
\]
So
\[\begin{split}
\|II\|_{L^1}
&\leq C\,\|g\|_{L^{q^*}}\,\sum_{i\leq |\log h|} 
 (\|A_i(.)\|_{L^q}+\|A_{i+1}(.)\|_{L^q}+\|A_{i-1}(.)\|_{L^q})\\
&\qquad\qquad\qquad \int_{B(0,1)}\int_{h}^{2^{-i}}
2^i\,r\,|z|\,\frac{dr}{r+h}\,\frac{dz}{|z|^{d-1}},  
\end{split}\]
where we used the localization in Fourier space of the $A_i$ and more
precisely the well known property
\[
\|A_i(.+\eta)-A_i(.)\|_{L^q}\leq C\,
2^i\,|\eta|\,(\|A_i(.)\|_{L^q}+\|A_{i+1}(.)\|_{L^q}+\|A_{i-1}(.)\|_{L^q}). 
\]

One then concludes that
\[
\|II\|_{L^1}\leq C\,\|g\|_{L^{q^*}}\,\sum_{i\leq |\log h|}
\|A_i\|_{L^q}=C\,\|g\|_{L^{q^*}}\,\|A\|_{B^0_{q,1}}.
\]
Combining the estimates on $I$ and $II$ in \eqref{decompQ1} gives
Lemma \ref{Q1}.

\medskip

For the proof of Lemma \ref{Q1l2}, it is enough to observe that in the
case $q=2$ by Cauchy-Schwartz
\[
\sum_{i\leq |\log h|} \|A_i\|_{L^2}\leq |\log
h|^{1/2}\,\left(\sum_{i} \|A_i\|_{L^2}^2\right)=|\log
h|^{1/2}\,\|A\|_{L^2}. 
\]
\subsection{Control on $Q_2$ for $A \in L^2$}
As for usual Calderon-Zygmund theory, the optimal
 bound on $Q_2$ is obtained in a $L^2$ setting namely
\begin{lemma}$\exists C>0$ s.t. for any $g$ and $A$ smooth with
  compact support
\[
\|Q_2(A,g)\|_{L^2}\leq C\,\|g\|_{L^\infty}\,\|A\|_{L^2}.
\]\label{Q2}
\end{lemma}

To prove this, first bound
\[
\begin{split}
|Q_2(A,g)|\leq \|g\|_{L^\infty}\,
\int_0^{1}\!\!\int_{S^{d-1}}\! \Bigg|&\int_{B(0,1)}\frac{r^{d}\, 
L(z,\omega)}{(r+h)^{d+1}}\,:\\
&\qquad\sum_{i\geq -\log_2 r}A_i(x+rz)\frac{dz}{|z|^{d-1}}\Bigg|\,d\omega\,\,dr.
\end{split}
\]
Hence
\[
\begin{split}
\|Q_2(A,g)\|_{L^2}\leq C\,R\,\|g\|_{L^\infty}^2+\frac{C}{R}\,
\int_0^{1}\!\!\int_{S^{d-1}}\! \Bigg\|&\int_{B(0,1)}\frac{r^{d}\, 
L(z,\omega)}{(r+h)^{d+1}}\,:\\
& \sum_{i\geq -\log_2
  r}A_i(.+rz)\frac{dz}{|z|^{d-1}}\Bigg\|_{L^2}^2
\,d\omega\,\,dr.
\end{split}
\]
Use Fourier transform and Plancherel equality on the last term to
bound it by
\[
\begin{split}
\sum_{\alpha,\beta=1}^3
\int_0^{1}&\frac{1}{r+h}\int_{S^{d-1}}\!\int_{|\xi|\geq 1/r} |{\cal
  F}A_{\alpha\beta}(\xi)|^2\\
& \int_{B(0,1)\times B(0,1)} L_{\alpha\beta}
(z,\omega)\,L_{\alpha\beta}(z',\omega)\,e^{i\xi\cdot
  r(z-z')}\frac{dz}{|z|^{d-1}}
\,\frac{dz'}{|z'|^{d-1}}\,d\xi\,d\omega\,dr
\end{split}
\]
One only has to bound the multiplier
\[
m(\xi,\omega,r)=\int_{B(0,1)\times B(0,1)} L_{\alpha\beta}
(z,\omega)\,L_{\alpha\beta}(z',\omega)\,e^{i\xi\cdot
  r(z-z')}\frac{dz}{|z|^{d-1}}
\,\frac{dz'}{|z'|^{d-1}}.
\]
Define
\[
M(\xi,\omega,r,s)=\int_{S^{d-1}} L_{\alpha\beta}(su,\omega)
e^{irs\xi\cdot u}\,du,
\]
such that
\[
m(\xi,\omega,r)=\int_0^1\int_0^1 M(\xi,\omega,r,s)\,\bar
M(\xi,\omega,r,s')
\,ds\,ds'.
\]
Assuming for instance that $\xi$ is along the first
axis, by the regularity on $\psi$ and hence $L$ given by Lemma
\ref{teclemma} 
\[\begin{split}
|M(\xi,\omega,r,s)|&=\left|\int_{S^{d-1}} L_{\alpha\beta}(su,\omega)
e^{irs|\xi| u_1}\,du\right|\leq \frac{1}{r\,|\xi|}
\int_{S^{d-1}} |\nabla_z L_{\alpha\beta}(su,\omega)|\,du\\
&\leq \frac{C}{r\,s\,|\xi|}.
\end{split}
\]
As $M$ is also obviously bounded, one deduces that
\[
|M(\xi,\omega,r,s)|\leq \frac{C}{\sqrt{r\,s\,|\xi|}}.
\]
Introducing this in $m$ immediately gives
\[
m(\xi,\omega,r)\leq \frac{C}{r\,|\xi|}.
\]
Therefore eventually
\[\begin{split}
\|Q_2(A,g)\|_{L^2}
&\leq C\,R\,\|g\|_{L^\infty}^2+\frac{C}{R}\,\sum_{\alpha,\beta=1}^3
\!\int_{\R^d} |{\cal
  F}A_{\alpha\beta}(\xi)|^2 \int_{|\xi|^{-1}}^1
\frac{1}{r+h}\,\frac{1}{r\,|\xi|}\,dr\,d\xi\\
&\leq C\,R\,\|g\|_{L^\infty}^2+\frac{C}{R}\,\sum_{\alpha,\beta=1}^3
\!\int_{\R^d} |{\cal
  F}A_{\alpha\beta}(\xi)|^2\,d\xi\\
&\leq C\,\|g\|_{L^\infty}\,\|A\|_{L^2},   
\end{split}\]
by optimizing in $R$, which proves the lemma.
\subsection{Control on $Q$ for $A \in L^p$}
To get an optimal bound, one should now try to obtain weak-type
estimates on $Q_2$, showing for instance that it belongs to $L^1-weak$
if $A\in L^1$; and then use interpolation. Additionally, we would have
to use the bound given by Lemma \ref{Q1} with Besov spaces. 

However here we will be satisfied with any explicit rate, even if it
is not optimal. We hence completely avoid some (not negligible)
technical difficulties and obtain
instead
\begin{lemma}  $\forall 1< q< \infty$, 
$\exists C>0$ s.t. for any
  smooth $g$ and $A$ with compact support  
\[
\|Q(A,g)\|_{L^1+L^q}\leq C\,|\log h|^{1/\bar
  q}\,\|g\|_{L^\infty\cap L^2}\,\|A\|_{L^q},
\]
where $\bar q=\min(q,q^*)$ with $1/q^*+1/q=1$.
\label{Qcomplete}
\end{lemma}
\noindent{\bf Remark.} Note that thanks to our decomposition of
$\nabla a$, we only use Lemma \ref{Qcomplete} for $q\leq 2$.

\medskip

\noindent{\bf Proof of Lemma \ref{Qcomplete}.} Fix $g$ and consider
$Q(A,g)$ as a linear operator on $A$. 
The easy control for $r\leq h$, Lemmas \ref{Q1l2} and \ref{Q2} imply that this
operator is bounded from $L^2$ to $L^1+L^2$ 
with norm $C\,(\|g\|_{L^\infty}+|\log h|^{1/2}\,\|g\|_{L^2})$.

On the other hand, one has the easy estimate
\[
|\bar Q(A,g)|\leq C\,\|g\|_{L^\infty}\,
\int_0^1 \frac{1}{r+h}\int_{B(0,1)} |A(x+rz)|\,\frac{dz}{|z|^{d-1}}.
\]
Therefore for any $1\leq q\leq \infty$, $\bar Q$ is bounded on $L^q$ with norm
less than $C\,\|g\|_{L^\infty}\,|\log
h|$. By usual interpolation, one deduces the lemma.
\subsection{Conclusion on the proof of Prop. \ref{theprop}}
By subsection \ref{begproof} 
\[\begin{split}
&\int_{\R^{2d}} \nabla K_h(x-y)\,(a(x)-a(y))\,|g(x)-g(y)|^2\,dx\,dy\\
&\quad \leq C\,(\|\hbox{div} a\|_{L^\infty}+\|\nabla a\|_{L^p})\,\int_{\R^{2d}}
K_h(x-y)\,|g(x)-g(y)|^2\,dx\,dy\\
&\quad\ +
2\|g\|_{L^1\cap L^\infty}\,\|Q(A,g)\|_{L^1+L^\infty}
+2\,\|g^2\|_{L^{p^*}}\,\|Q(A,1)\|_{L^p}.
\end{split}
\]
Bound
directly $Q(A,g)$ by Lemma \ref{Qcomplete} and observe that $Q(A,1)$
is bounded on any $L^p$ with $1<p<\infty$. This completes the proof of
the proposition.

\medskip

\noindent{\bf Acknowledgement.} P-E Jabin is grateful to
P.A. Markowich for having pointed out this interesting problem.


\end{document}